\documentclass{llncs}
\usepackage{stmaryrd}
\usepackage{amsmath}
\usepackage{amssymb}
\usepackage{enumerate}
\begin{document}

\title{MEANING, CHOICE and ALGEBRAIC SEMANTICS of SIMILARITY BASED ROUGH SET THEORY}

\author{A. MANI}
\institute{Member, Calcutta Mathematical Society\\
9/1B, Jatin Bagchi Road\\
Kolkata(Calcutta)-700029, India\\
\email{$a.mani@member.ams.org$}\\
Home Page: \url{http://amani.topcities.com}}

\titlerunning{Meaning, Choice and Algebraic Semantics}

\maketitle

\begin{abstract}
Both algebraic and computational approaches for dealing with similarity spaces are well known in generalized rough set theory. However, these studies may be said to have been confined to particular perspectives of distinguishability in the context. In this research, the essence of an algebraic semantics that can deal with all possible concepts of distinguishability over similarity spaces is progressed. Key to this is the addition of choice-related operations to the semantics that have connections with modal logics as well. Among these we focus on a semantics that is based on \emph{local clear distinguishability} over similarity spaces.  
\end{abstract}
\section{Introduction}
By a \emph{Tolerance Approximation Space} (TAS), we mean a pair of the form $S\,=\,\left\langle \underline{S},\,T\right\rangle$, with $\underline{S}$ being a set and $T$ a tolerance relation over it - these are also known as similarity or tolerance spaces. Some references for extension of classical rough set theory to TAS are \cite{KM} \cite{CG98}, and \cite{KB98}. The type of granules used in these theories is summarized below.

An approach (\cite{KM}) has been to define a new equivalence $\theta_{0}$ on  $S$ via $(x,\,y)\,\in\,\theta_{0}$ if and only if $dom_{T}(x)\,=\,dom_{T}(y)$  with  $dom_{T}(z)\,=\,\cap\{[x]_{T}\,:\,z\,\in\,[x]_{T}\}$. This is essentially an unduly cautious 'clear perspective' approach.

A somewhat natural generalization of the approximation space semantics using $T$-related sets (or tolerance sets) can be described from the point of view of generalized covers (see \cite{IM}). This includes the approach of defining  $[x]_{T}\,=\,\{y\,;\,(x,\,y)\,\in\,T\}$ and the lower and upper approximation of a set $A$ as, $A^{l}\,=\,\bigcup\{[x]_{T}\,;\,[x]_{T}\,\subseteq\,A\}$ and $A^{u}\,=\,\bigcup\{[x]_{T}\,;\,[x]_{T}\,\cap\,A\,\neq\,\emptyset ,\,x\,\in\,A\}.$
A \emph{bited modification} proposed in \cite{SW}, valid for many definable concepts of granules, consists in defining a \emph{bited upper approximation}. Algebraic semantics of the same has been considered by the present author in \cite{AM105}. It is also shown that a full representation theorem is not always possible for the semantics.    

In \cite{CG98}, the approximations $A^{l*}\,=\,\{ x\,;\,(\exists{y})\,(x,\,y)\,\in\,T,\,[y]_{T}\,\subseteq\,A\}$ and 
$A^{u*}\,=\,\{x\,;\,(\forall{y})\,((x,\,y)\,\in\,T\,\longrightarrow\,[y]_{T}\,\cap\,A\,\neq\,\emptyset)\}.$ are introduced. It can be shown that, for any subset $A$, $A^{l}\subseteq\,A^{l*}\,\subseteq\,A\,\subseteq\,A^{u*}\,\subseteq\,A^{u}$. 

In the BZ andQuasi-BZ algebraic semantics (\cite{CC}), the lower and upper rough operators are generated by a preclusivity operator and the complementation relation on the power set of the approximation space, or on a collection of sets under suitable constraints in a more abstract setting. Semantically the BZ-algebra and variants do not capture all the possible ways of arriving at concepts of distinguishability over similarity spaces. Whereas the quasi-BZ lattice does not encompass a paradigm shift relative the BZ-algebra, the BZMV variants are designed to capture fuzzy aspects. 

In subjective terms, reducts are minimal sets of attributes that preserve the quality of classification. An important problem is in getting good scalable algorithms for the computation of the different types of reducts (or supersets that are close to them) (see ~\cite{KC}). These depend on the concept of granules used. For TAS, most of the above concepts of granules and approximation often lead to computational difficulties. Though our theory can be expected to improve the computational situation, the main motivation for the present work is centred around our new concept of \emph{local clear distinguishability}.

\begin{definition}
Let $P\,=\,\left\langle \underline{P},\,<\,\right\rangle$ be a partially ordered set and if $A$ is any subset of $P$, let its lower and upper cone be $L(A)\,=\,\{x\,;\,(\forall a\,\in\,A)\,x\,\leq\,a\}$ and $U(A)\,=\,\{x\,;\,(\forall a\,\in\,A)\,a\,\leq\,x\}$ respectively. A function $\lambda\,:\,\wp(P)\,\mapsto\,P$ will be said to be \emph{lattice-coherent} with $<$ if and only if 
the condition $a\,\leq\,b$ then $\lambda (L(a,\,b))\,=\,a$ and $\lambda (U(a,\,b))\,=\,b$.

By a \emph{choice function} $\chi$ on a set $S$, we mean a function $\chi\,:\wp(S)\,\longmapsto\,S$, which is such that $(\forall x\,\in\,S)\,\chi (\{x\})\,=\,x$ and $(\forall A\,\in\,\wp(S))\,\chi (A)\,\in\,A$.
\end{definition} 

\section{Philosophical Basis}

In this research we deal with choice-based rough granulation as opposed to basing choice forms over rough granulation. Choice is needed to specify the admissible granulation in the theory. A cautious way of saying that a set of things are essentially indistinguishable is to say that they are mutually indistinguishable. If we are to stick to this principle, then given a subset $A$ of a TAS $S$, the essentially indistinguishable subsets of $A$ are the intersections of the blocks of the relation $T$ with $A$. In course of constructing an \emph{upper approximation} of $A$, we can form all the unions of disjoint blocks (as opposed to \emph{union of all blocks}) that intersect with $A$. The operation enshrines a natural concept of clarity in the process of categorisation. The other option fails to do so. This preference can be viewed as a local equivalence-based perspective in the context. It is implicit that 'at a local level' the concept of distinguishability is based on disjoint categories. The loss of uniqueness in the construction can be dealt with through suitable choice functions that are latent in the context. 

By the \textbf{local clear distinguishability} principle (LCP), we mean the requirement that definite objects generated by approximation context initiators should be made up of nonintersecting granules. The most appropriate domain of discourse for general rough set theory and this concept should avoid ZF sets, but for simplicity's sake we use the latter.  A context initiator may be the subset under consideration or a variant thereof depending on the application context. The \emph{local} part is, because we restrict to granules generated from these objects. For example, if we perceive a subset $A\,=\,\{a_{1},\,a_{2},\,\ldots\,a_{n}\}$ of a TAS $S$ in the given order, and use a 'First In First Out' (FIFO) principle in the generation of a maximal set of disjoint blocks contained within $A$, then the final outcome is determined by the 'order of perception' and the 'generation principle'. I will consider a perfect technical formulation of the concept in a separate paper as the fine details are not explicitly required in our theory.

In learning theory, it can be interesting to know the structure of contextual knowledge that is guided by the local clear distinguishability principle. The latter in essence is a strategy for forming a clear concept of context-dependent knowledge. This is in line with Pawlak's concept of knowledge in classical rough set theory, where if $\underline{S}$ is a set of attributes, then sets of the form $A^{l}$ and $A^{u}$ represent clear and definite concepts. If $Q$ is another stronger equivalence on $\underline{S}$, then the state of the knowledge encoded by $\left\langle \underline{S},\,Q \right\rangle $ is a \emph{refinement} of that of $S\,=\,\left\langle\underline{S},\,P\right\rangle $. 

The entire set $S$ will not be a union of disjoint blocks in general and we will be able to find maximal collections of mutually disjoint blocks. Deciding on what ought to be the upper approximation of $S$ and the following are necessary. If a set that seems to be an upper approximation of another given set $A$ (say) on the basis of nonempty intersection of blocks with $A$, fails to satisfy the local clarity paradigm on application of a strategy similar to that used for the lower approximation, then we can 
\begin{itemize}
\item {Take the upper approximation of $A$ to be undefined}
\item {Relax the local clarity paradigm and take the upper approximation as $S$}
\item {Relax the inclusion of the set in its upper approximation by way of selecting a union of a set of disjoint blocks from the set of all blocks that intersect the set $A$ in question.}
\item {Take the upper approximation to be one of the unions of maximal collections of mutually disjoint blocks}
\item {Relax the local clarity paradigm and take the relevant upper approximation to be the same as the union of blocks that intersect the set $A$ in question.}
\end{itemize}

In TAS, it can be suspected that approximations of the above kind actually improve the information content of possible semantics to the point that we have good representation theorems as well. In this paper we develop an elegant semantics with the the first option for the concept of an upper approximation (for the other options see \cite{AM99}) and show that this is indeed the case. Interestingly the restriction that \emph{if a set is a union of disjoint blocks, then it ought to be exact (or crisp)} in conjunction with the above choice scenario turns out to be the basis of a nonmonotone variation of the theory (\cite{AM99}). In all this, the nature of the choice process consists in selecting particular subcollections of mutually disjoint blocks from sets of such collections. 

It can be argued that the approximation contexts generated by TAS should involve more than the concepts of lower and upper approximation and perhaps a gradation of the concept. Our semantic approach actually supports this and in fact we think that anything less than four must involve loss of information. Also any reasonable gradation must necessarily depend on the topology on the TAS or a variant of the notion. But we will not explicitly refer to such structures in this paper. 

Rough approximation in its general forms is distinct from \emph{approximation} in its ability to include clear means of categorisation in the approximation. So in comparing different types of rough approximation, a simple concept of \emph{fineness} of the approximation can never be a suitable criteria for differentiating between approximation methodologies. The lower approximation of a set using the equivalence $\theta_{0}$ mentioned in the introduction will be \emph{very close} to the set and so will the upper approximation be. 'Closeness to the set in question' is not a sufficient criteria for deciding among the concepts of approximation as it often ends up violating a variety of context-dependent coherence criteria. Such insufficient criterias have been used often in the literature. 

\section{Essential $\lambda$-Rough Partial Algebras}

Let $S\,=\,\left\langle \underline{S},\,T\right\rangle$ be a TAS, $A\,\subseteq\,\underline{S}$ be an arbitrary subset of it and let $\mathbb{S}$ be the collection of all blocks of $T$. We can endow $\wp(\mathbb{S})$ with the $\prec$ order. The $\prec$ order being defined via, if $E,\,B\,\in\,\wp(\mathbb{S})$ then $E\,\prec\,B$ if and only if $E\,\subseteq\,B$ and $E$ is a subcollection of disjoint blocks. 
\begin{description}
\item[Lower Relativisation] {Form the collection $\mathbb{S}(A)$ of all blocks included in $A$}
\item[Lower Clarification-1,2] {Form the collection $\mathbb{LS}(A)$ of subcollections of mutually disjoint elements in $\mathbb{S}(A)$, order these by inclusion and determine the collection of maximal elements $\mathbb{LS}_{M}(A)$.}
\item[Choice] {We will assume that we have a choice function $\lambda:\,\wp(\wp(\mathbb{S}))\,\mapsto\,\wp(\mathbb{S})$ that is lattice-coherent with the $\prec$ order on the collection $\wp(\mathbb{S})$. }
\item[Lower Choice] {$\bigcup \lambda(\mathbb{LS}_{M}(A))$ will be called the \emph{0-lower approximation} of $A$. It will be abbreviated to $A^{l0}$.}
\item[Primitive Lower Choice] {$\lambda(\mathbb{LS}_{M}(A))$ will be called the \emph{primitive lower approximation} of $A$}
\item[Lateral Lower Choice]{$\bigcup \mathbb{S}(A)$ will be called the \emph{lateral lower approximation} of $A$ and will be denoted by $A^{\breve{l}}$}
\item[Upper Relativisation] {Form the collection $\mathbb{S}_{u}(A)$ of all blocks that intersect with $A$. }
\item[Upper Clarification-1, 2] {Let $\mathbb{US}_{m}(A))$ be the set of minimal elements in the set of subcollections of mutually disjoint blocks in $\mathbb{S}_{u}(A)$ each of whose unions contains $A$.}
\item[Upper Choice] {$\bigcup \lambda(\mathbb{US}_{m}(A))$ will be called the \emph{0-upper approximation} of $A$. It will be abbreviated to $A^{u0}$. If $\mathbb{S}_{u}(A)$ is empty, then take $A^{u0}$ to be undefined.}
\item[Primitive Upper Choice] {$\lambda(\mathbb{US}_{m}(A))$ will be called the \emph{primitive upper approximation} of $A$}
\item[Lateral Upper Choice] {$\bigcup \mathbb{S}_{u}(A)$ will be called the \emph{lateral upper approximation} of $A$ and will be denoted by $A^{\breve{u}}$}
\end{description}

\begin{theorem}
All of the above approximations are all well-defined and satisfy the following properties: 
\begin{enumerate}[(a)]
\item {For any subset $A$, $(A^{l0})^{l0}\,=\,A^{l0}\,\subseteq\,A^{\breve{l}}$}
\item {For any subset $A$, $A^{l0}\,\subseteq\,(A^{l0})^{u0}$}
\item {For any subset $A$, $(A^{u0})^{l0}\,\stackrel{w}{=}\,A^{u0}\,\stackrel{w}{=}\,(A^{u0})^{u0}\,\subseteq\,A^{\breve{u}}$ ; For terms $p,\,q$, $p\stackrel{w}{=}q$ iff $(\forall x\,\in\,dom(p)\,\cap\,dom(q)) p(x)\,=\,q(x)$ (of course w.r.t an interpretation)}
\item {$(A\,\subseteq\,B\,\longrightarrow\,A^{l0}\,\subseteq\,B^{l0})$}
\item {$(A\,\subseteq\,B,\,A\,\subseteq\,A^{u0}\,B\,\subseteq\,B^{u0}\longrightarrow\,A^{u0}\,\subseteq\,B^{u0})$}
\item {$(A\,\subseteq\,B\,\longrightarrow\,A^{\breve{l}}\,\subseteq\,B^{\breve{l}},\,A^{\breve{u}}\,\subseteq\,B^{\breve{u}})$}
\item {If $A$ is a subset of a TAS $S$ and $A^{l0}\,=\,A\,=\,A^{\breve{l}}$, then $A$ is necessarily a union of disjoint blocks.}
\item {If $A^{u0}$ exists, then
$A^{l0}\,\subseteq\,A^{l}\subseteq\,A^{l*}\,\subseteq\,A^{l\theta}\,\subseteq\,A\,\subseteq\,A^{u\theta}\,\subseteq\,A^{u0}\,\subseteq\,A^{u*}$
else, $A^{l0}\,\subseteq\,A^{l}\subseteq\,A^{l*}\,\subseteq\,A^{l\theta}\,\subseteq\,A\,\subseteq\,A^{u\theta}\,\subseteq\,A^{u*}$}
\item {If $A$ is a subset of $S$ that is also a block of the tolerance, then $A^{l0}\,=\,A\,=\,A^{\breve{l}}$, but it can happen that $A^{u0}\,\neq\,A^{l0}$ and $A^{\breve{u}}\,\neq\,A^{u0}$ }
\end{enumerate}
\end{theorem}
\begin{proof}
In general if we apply the 0-upper approximation construction to a set of the form $A^{l0}$, then we will get a larger set. The other parts can be verified by direct set-theoretic arguments. For the last two claims note that, $A^{l*}\,=\,\{ x\,;\,(\exists{y})\,(x,\,y)\,\in\,T,\,[y]_{T}\,\subseteq\,A\}$\\

$A^{u*}\,=\,\{x\,;\,(\forall{y})\,((x,\,y)\,\in\,T\,\longrightarrow\,[y]_{T}\,\cap\,A\,\neq\,\emptyset)\}$ and that a block cannot contain any other blocks. \qed
\end{proof}

Note that the property $A^{l0}\,\subseteq\,(A^{l0})^{u0}$ also happens in esoteric rough set theory [\cite{AM24}]. When we redefine situations in which the 0-upper approximations are undefined with those values set to $S$, then the whole of the behaviour (as far as the two approximations $l0$ and the new $u0$ are concerned) resembles that of rough set theory over partial approximation spaces.

Lateral approximations do not encompass discernibility at the local level. So we do not think that they constitute a reasonable rough concept by themselves. The other non modal approximations indicated in the introduction are still distinct from these. They can be obtained through suitable modifications.

Parts-(g) and (h) of the above therem have a deep role to play in deciding on the direction of possible representation theorems. The propositions ensure that we can identify unions of disjoint blocks using the approximation operators.

\begin{theorem}
If we define the operations $\sim,\,\minuso$ over the power set $\wp (S)$ via (the latter being a partial operation that is defined only when $A^{u0}$ is)
\[\sim A\,=\,S\,\setminus\,A^{\breve{u}} \;\;\; \minuso A\,=\,S\,\setminus\,A^{u0},\] then it is necessary that $A\,\subseteq\,\sim\sim A$, but in general $A\,\nsubseteq\,\minuso \minuso A$, even when the right hand side is defined. 
\end{theorem}
\begin{proof}
Suppose the contrary, $x\,\in\,A$ and $x\,\notin\,\sim\sim A$, then $x\,\in\,A^{\breve{u}}$. This means $x\,\notin\,S\,\setminus\,A^{\breve{u}}$. As $x\,\in\,(S\,\setminus\,A^{\breve{u}})^{\breve{u}}$ (by assumption), there must exist a block $F$ such that $x\,\in\,F$ and $F\,\cap\,(S\,\setminus\,A^{\breve{u}})\,\neq\,\emptyset$. Since $x\,\in\,F\,\cap\,A$, so $F\,\subseteq\,A^{\breve{u}}$. This contradiction implies that the original assumption must be false and therefore $A\,\subseteq\,\sim\sim A$. Counterexamples for the second part are easy. \qed
\end{proof}

\begin{proposition}
For any $A,\,B\,\subseteq\,S$, let $A\,\curlyvee\,B\,=\,(A^{u0}\,\cup\,B^{u0})^{u0}$ (if defined) and $A\,\curlywedge\,B\,=\,(A^{u0}\,\cap\,B^{u0})^{l0}$ (if defined) then the following holds for the partial operations:
\begin{enumerate}[(a)]
\item {$A\,\curlyvee\,A\,\stackrel{w}{=}\,A^{u0 u0}$}
\item {$(A\,\curlyvee\,B)^{l0}\,=\,A\,\curlyvee\,B$}
\item {If $\mathbb{S}$ is the set of two element subsets of $S$ that are not included in any block, then  $A\,\in\,\mathbb{S}\,\longrightarrow\,A^{l0}\,=\,\emptyset\,A^{\breve{l}}$ and $A^{u0}$ is undefined.}
\item {If $\mathbb{S}$ is the set of two element subsets of $S$ that are not included in any block, then $A^{\breve{u}}$ is a union of at least two blocks.}
\item {$A\,\curlywedge\,A\,\stackrel{w}{=}\,A^{u0}\,; \;\;A\,\curlywedge\,B\,\stackrel{w}{=}\,B\,\curlywedge\,A$}
\item {$(A\,\subseteq\,B\,\subseteq\,B^{u0}\longrightarrow\,A\,\curlyvee\,B\,=\,B^{u0})$}
\end{enumerate} 
\end{proposition}
\begin{proof}
The proof is by direct arguments. Note that if a subset is not contained in any block, then it must contain a two element subset that has the same property. This motivates the third and fourth claims. \qed  
\end{proof}

\begin{definition}
A \emph{pre-essential $\lambda$-rough partial algebra} will be an algebra of the form $\Xi(S)\,=\,\left\langle\underline{\wp(S)|\sigma},\,\leq,\,\curlyvee,\,\curlywedge,\,\sqcup,\, \ovee,\,\owedge,\,L_{0},\,U_{0},\,\breve{L},\,\breve{U},\sim,\,\minuso,\,[\emptyset],\,[S] \right\rangle $ that has been constructed as follows from a TAS $S$:
\begin{itemize}
\item {For any set $A\,\in\wp(S)$, if $A^{u0}$ is defined let $\upsilon (A)\,=\,(A^{l0},\,A^{u0},\,A^{\breve{l}},\,A^{\breve{u}})$, else let $\upsilon (A)\,=\,(A^{l0},\,A^{\breve{l}},\,A^{\breve{u}})$}
\item {Let $(A,\,B)\,\in\,\sigma$ if and only if $\upsilon(A)\,=\,\upsilon(B)$ }
\item {Then form the quotient $\wp(S) | \sigma$}
\item {Define $L_{0}([A])\,=\,[A^{l0}]$, $U_{0}([A])\,=\,[A^{u0}]$ if defined}
\item {On the quotient, let $[A]\,\leq\,[B]$ if and only if $A^{l0}\,\subseteq\,B^{l0}$ and $A^{u0}\,\subseteq\,B^{u0}$ (if defined) and $A^{\breve{l}}\,\subseteq\,B^{\breve{l}}$ and $A^{\breve{u}}\,\subseteq\,B^{\breve{u}}$. We will denote the strict version of the inequality by $\lneq$ }
\item {Define $[A]\,\ovee\,[B]\,\stackrel{\text {\tiny def}}{=}\,[A^{u0}\,\cup\,B^{u0}]$ if defined}
\item {Define $[A]\,\owedge\,[B]\,\stackrel{\text {\tiny def}}{=}\,[A^{u0}\,\cap\,B^{u0}]$ if defined}
\item {Define $[A]\,\curlyvee\,[B]\,\stackrel{\text {\tiny def}}{=}\,U_{0}([A]\,\ovee\,[B])$ if defined}
\item {Define $[A]\,\curlywedge\,[B]\,\stackrel{\text {\tiny def}}{=}\,L_{0}([A]\,\owedge\,[B])$ if defined}
\item {Define $\breve{U}([A])\,\stackrel{\text {\tiny def}}{=}\,[A^{\breve{u}}]$,  $\breve{L}([A])\,=\,[A^{\breve{l}}]$}
\item {Define $[A]\,\sqcup\,[B]\,\stackrel{\text {\tiny def}}{=}\,[A\,\cup\,B]$}
\item {Define $[A]\,\sqcap\,[B]\,\stackrel{\text {\tiny def}}{=}\,[A\,\cap\,B]$}
\item {Define $\sim[A]\,\stackrel{\text {\tiny def}}{=}\,[S\,\setminus\,A^{\breve{u}}]$}
\item {Define $\minuso[A]\,\stackrel{\text {\tiny def}}{=}\,[S\,\setminus\,A^{u0}]$ if defined}
\end{itemize}
\end{definition}

\begin{theorem}
All of the fundamental and derived operations of a pre-essential $\lambda$-rough partial algebra are well defined and all of the following hold:
\begin{enumerate}[(a)]
\item {If $x$ is a class corresponding to a union of disjoint blocks then $L_{0}x\,=\,\breve{L}x\,=\,x$ and conversely.}
\item {If $x\,=\,U_{0}x$, then $x$ is a class corresponding to a union of disjoint blocks, but the converse need not hold in general.}
\item {If $x$ is a class generated by a single block, then $L_{0}x\,=\,U_{0}x\,=\,x\,=\,\breve{L}x $ }
\item {If for a class $x$, $(\forall y) (y\, \lneq\,x\,\longrightarrow\,y\,\neq\,L_{0}(y)\,\neq\,U_{0}(y))$ and $L_{0}(x)\,=\,U_{0}(x)\,=\,x$, then $x$ is the class corresponding to a single block and conversely }
\item {If for a class $x$, $(\forall y) (y\, \lneq\,x\,\longrightarrow\,L_{0}(y)\,\lneq\,L_{0}(x))$ and $L_{0}(x)\,=\,U_{0}(x)\,=\,x$, then $x$ is the class corresponding to a single block and conversely}
\item {If for a class $x$ that does not correspond to that of a single block, ($L_{0}(x)\,=\,[\emptyset]$ or $L_{0}(x)$ corresponds to a single block) and $U_{0}(x)$ is undefined, then $x$ is a class that corresponds to a set that contains a two element set that is not in any block of $T$.}
\item {$(U_{0}(x)\,=\,U_{0}(x)\,\longrightarrow\,\minuso (x)\,\leq\,\minuso (L_{0}(x)))$}
\end{enumerate} 
\end{theorem}

\begin{proof}
A partial operation $f\:X^{n}\,\mapsto\,X$ is well defined if at each point it is uniquely defined or not non-uniquely defined at all. Most of the proof is included in the proof of the theorem for essential $\lambda$-rough partial algebras below. \qed
\end{proof}

\begin{definition}
In the light of the above theorem, we introduce the following derived operations and predicates on a pre-essential $\lambda$-rough algebra $\Xi(S)$,
\begin{itemize}
\item {For any $x$, if $x$ is the class of a single block, then let $s(x)\,=\,x$, else $s(x)\,=\,\emptyset$.}
\item {For any $x$, if $x$ is the class of a 2-element subset that is not in any block, then let $t(x)\,=\,x$, else $t(x)\,=\,\emptyset$.}
\item {$IU(x)$ if and only if $U_{0}(x)\,=\,U_{0}(x)$. Note that $U_{0}$ is a partial operation. Further we will write $IU(a, b,..)$ for $IU(a),\, IU(b),\ldots$}
\item {$IN(x)$ if and only if $\minuso x\,=\,\minuso x$.}
\end{itemize}
The algebra formed by adjoining the additional operations and predicates ($\lneq,\, s,\,t,\,IU,\,IN $) to $\Xi(S)$ will be termed an \emph{essential $\lambda$-rough partial algebra} and denoted by $\maltese(S)$.  
\end{definition}

\begin{theorem}
All of the following hold in an essential $\lambda$-rough partial algebra $\maltese (S)$:
\begin{enumerate}[(a)]
\item {$x\,\curlyvee\,y\,\stackrel{w}{=}\,y\,\curlyvee\,x$ ; $x\,\curlywedge\,y\,\stackrel{w}{=}\,y\,\curlywedge\,x$ ; $x\,\curlyvee\,x\,\stackrel{w}{=}\,U_{0} U_{0}(x)$ ; $x\,\curlywedge\,x\,\stackrel{w}{=}\,U_{0}(x)$}
\item {$L_{0}(x)\,\leq\,\breve{L}(x)\,\leq\,x\,\breve{U}(x)$ ; $(IU(x)\,\longrightarrow\,x\,\leq\,U_{0}(x))$}
\item {$L_{0}L_{0}(x)\,=\,L_{0}(x)$ ; $IU(x)\,\longrightarrow\,U_{0}(x)\,\leq\,U_{0}U_{0}(x)$}
\item {$\breve{L}L_{0}(x)\,=\,L_{0}(x)$ ; $L_{0}\breve{L}(x)\,\leq\,\breve{L}(x)$}
\item {$\breve{L}\breve{L}(x)\,=\,\breve{L}(x)$ ; $(IU(x)\,\longrightarrow\,L_{0}U_{0}(x)\,=\,U_{0}(x))$}
\item {$L_{0}(x)\,\leq\,U_{0}L_{0}(x)$ ; $\breve{U}(x)\,\leq\,\breve{U}\breve{U}(x)$}
\item {$(IU(x)\,\longrightarrow\,x\,\leq\,U_{0}(x)\,\leq\,\breve{U}(x)\,\leq\,\breve{U}U_{0}(x)\,\leq\,\breve{U}\breve{U}(x))$}
\item {$(x\,\leq\,y\,\longrightarrow\,L_{0}(x)\,\leq\,L_{0}(y),\,\breve{U}(x)\,\leq\,\breve{U}(y),\,\breve{L}(x)\,\leq\,\breve{L}(y))$}
\item {$(x\,\leq\,y,\,IU(x)\,\longrightarrow\,U_{0}(x)\,\leq\,U_{0}(y))$}
\item {$(x\,\leq\,y,\,IU(x,\,y)\,\longrightarrow\,x\,\curlywedge\,y\,=\,U_{0}(x)\,=\,x\,\owedge\,y,\,x\,\curlyvee\,y\,=\,U_{0}(y)\,=\,x\,\ovee\,y)$}
\item {$(IU(x,\,y,\,a,\,b),\,x\,\leq\,y,\,a,\,\leq\,b\,\longrightarrow\,x\,\curlywedge\,a\,\leq\,y\,\curlywedge\,b)$}
\item {$(IU(x,\,y,\,a,\,b,\,x\,\ovee\,a,\,y\,\ovee\,b),\,x\,\leq\,y,\,a,\,\leq\,b\,\longrightarrow\,x\,\curlyvee\,a\,\leq\,y\,\curlyvee\,b)$}
\item {$t(x)\,=\,x$ if and only if $\neg(IU(x)),\,s(\breve{U}(x))\,=\,0,\,L_{0}(x)\,\lneq\,x,\,(\forall{y})(y\,\lneq\,x\,\longrightarrow\,ty\,=\,0)$ and $(0\,\lneq\,a,\,b,\,c\,\lneq\,x\,\longrightarrow\,a\,=\,b\; \mathrm{or} \;b\,=\,c\; \mathrm{or} \;c\,=\,a)$}
\item {$(IU(x)\,\longrightarrow\,\sim \minuso x\,\leq\,\minuso \sim x)$}
\item {$\sim x \,\leq\,\sim L_{0}(x)$ ; $\sim 0\,=\,1$ ; $\sim 1 \,=\,0$}
\item {$x\,\leq\,\sim \sim x$ ; $(IU(x)\,\longrightarrow\,\sim U_{0} (x)\,\leq\,\sim x)$}
\item {$\sim \breve{U}(x)\,\leq\,\sim x $ ; $(IU(x)\,\longrightarrow\,\sim \breve{U} (x)\,\leq\, \sim U_{0}(x)$}
\item {$\sim x \,\leq\,\sim \breve{L}(x)\,\leq\,\sim L_{0}(x)$ ; $(IN(x)\,\longrightarrow\,\minuso x\,\leq\,\minuso L_{0}(x))$}
\item {$(IU(\sim x),\,IU(x)\,\longrightarrow\,\sim \minuso x \,\leq\, \minuso \sim x)$ ; $\minuso 0\,=\,1$ ; $\neg IN(1)$}
\item {$(IU(x)\,\longrightarrow\,\minuso x\,\leq\,\minuso L_{0}(x))$ ; $(IN(x)\,\longrightarrow\,\minuso \breve{U}(x)\,\leq\,\minuso x,\,\minuso U_{0}(x)\,=\,\minuso x)$}
\item {$\forall{y}(y\,\lneq\,x\,\longrightarrow\,L_{0}(y)\,\lneq\,L_{0}(x))$ and $L_{0}(x)\,=\,U_{0}(x)\,=\,x$ if and only if $s(x)\,=\,x$, else $s(x)\,=\,0$}
\item {$x\,\sqcup\,y\,=\,y\,\sqcup\,x$ ; $x\,\sqcup\,x\,=\,x$ ; $x\,\ovee\,y\,\stackrel{w}{=}\,U_{0}(x)\,\sqcup\,U_{0}(y)$}
\item{$(x\,\leq\,y\,\longrightarrow\,x,\,y\,\leq\,x\,\sqcup\,y)$ ; $\breve(U)(x)\,\sqcup\,\breve{U}(y)\,\leq\,\breve{U}(x\,\sqcup\,y)$}
\end{enumerate}
\end{theorem}

\begin{proof}
Since an essential $\lambda$-rough partial algebra is a concrete object, we will assume $x\,=\,[A],\,y\,=\,[B]$ for some suitable subsets $A,\,B$.
\begin{enumerate}[(a)]
\item {If either side is defined, then the other is defined and is equal to $U_{0}[A^{u0}\,\cup\,B^{u0}]\,=\,[(A^{u0}\,\cup\,B^{u0})^{u0}]$. The next part follows in the same way.
If defined, $x\,\curlyvee\,x\,=\,[(A^{u0}\,\cup\,A^{u0})^{u0}]\,=\,[(A^{u0})^{u0}]\,=\,U_{0}U_{0}x $.  The last part follows in the same way.}
\item {For any element $A$ of a class $x$, $L_{0}(x)$ is the class of $A^{l0}$, while $\breve{L}(x)$ is the class of $A^{\breve{l}}$. $L_{0}(x)\,\leq\,\breve{L}(x)$ follows from the previously proved properties of the approximations as applied to the sets $A^{l0}$ and $A^{\breve{l}}$. $\breve{L}(x)\,\leq\,x\,\breve{U}(x)$ follows in the same way. 

In the second part $(IU(x)$ means that for each of the elements of $x$, the 0-upper approximation exists and they must all be equal to each other. Now if $A\,\in\,x$, then $A^{u0}\,\in\,U_{0}(x)$ and if $B\,\in\,U_{0}(x)$, then $(A^{u0})^{*}\,\subseteq\,B^{*}$ , where $*$ is any of the four approximations. So we have $IU(x)\,\longrightarrow\,x\,\leq\,U_{0}(x))$.}
\item {If $A\,\in\,L_{0}L_{0}(x)$ then it must be of the form $B^{l0}$ for some $B\,\in\,L_{0}(x)$, as for any $B\,\in\,L_{0}(x)$ we must have $A^{l0}\,=\,B^{l0 l0 l0}\,=\, B^{l0}$. The equality $ L_{0} L_{0} (x)\,=\,L_{0}(x)$ follows as a consequence. If $U_{0}(x)$ is defined, then $U_{0}U_{0}(x)$ will also be defined, for any $A\,\in\,x$, it is possible that $A^{u0}\,\subset\,A^{u0u0}$, so   $IU(x)\,\longrightarrow\,U_{0}(x)\,\leq\,U_{0}U_{0}(x)$}
\item {If $A\,\,\in\,L_{0}(x)$, then it is already a union of disjoint blocks and $A^{\breve{l}}$ will be equal to it, so $A\,\in\,\breve{L}L_{0}(x)$. As the converse inclusion is trivial so $\breve{L}L_{0}(x)\,=\,L_{0}(x)$. The second part follows by a similar argument. }
\item {If $A\,\in\,\breve{L}(x)$, then $A^{\breve{l}}\,\in\,\breve{L}\breve{L}(x)$ as $A^{\breve{l}\,=\,A^{\breve{l}\breve{l}}}$. Now if $B\,\in\,\breve{L}\breve{L}(x)$, then there is a $C\,\in\,\breve{L}(x)$ such that $B^{\breve{l}}\,=\,C^{\breve{l}\breve{l}}$, $B^{l0}\,=\,C^{\breve{l}{l0}}$ and so on for the other approximations. But the first equality means that $B^{\breve{l}}\,=\,C^{\breve{l}}$. So $\breve{L}\breve{L}(x)\,=\,\breve{L}(x)$.

If defined, elements of $U_{0}(x)$ are themselves identical unions of disjoint blocks, each of whose 0-lower approximations are the same. So the conclusion $L_{0}U_{0}(x)\,=\,U_{0}(x)$ follows} 
\item {Elements of $L_{0}(x)$ are identical unions of disjoint blocks whose 0-upper approximations are equal, this ensures that $U_{0}L_{0}(x)$ is defined. So $L_{0}(x)\,\leq\,U_{0}L_{0}(x)$.

$\breve{U}(x)\,\leq\,\breve{U}\breve{U}(x)$ is a special case of $x\,\leq\,\breve{U}\breve{U}(x)$}
\item {In $(IU(x)\,\longrightarrow\,x\,\leq\,U_{0}(x)\,\leq\,\breve{U}(x)\,\leq\,\breve{U}U_{0}(x)\,\leq\,\breve{U}\breve{U}(x))$, the premise is necessary for ensuring that $U_{0}(x)$ is defined.}
\item {If $x\,=\,[A]$ and $y\,=\,[B]$ for some $A,\,B$, then $x\,\leq\,y$ implies that $A^{l0}\,\subseteq\,B^{l0},\,A^{u0}\,\subseteq\,B^{u0},\,A^{\breve{l}}\,\subseteq\,B^{\breve{l}},A^{\breve{u}}\,\subseteq\,B^{\breve{u}}$. So the conclusion   $(x\,\leq\,y\,\longrightarrow\,L_{0}(x)\,\leq\,L_{0}(y),\,\breve{U}(x)\,\leq\,\breve{U}(y),\,\breve{L}(x)\,\leq\,\breve{L}(y))$ follows}
\item {$(x\,\leq\,y\,IU(x)\,\longrightarrow\,U_{0}(x)\,\leq\,U_{0}(y))$ is not part of the above statement, but can be easily verified.}

\item {If $x\,\leq\,y,\,IU(x,\,y$, then $Wx\,\leq\,Wy$, where $W$ is any of the four approximation operators. So then if $A\,\in\,x$ and $B\,\in\,B$, $x\,\ovee\,y\,=\,[A^{u0}\,\cup\,B^{u0}]\,=\,[B^{u0}]\,=\,U_{0}(y)$. 

In a similar way the rest of $(x\,\leq\,y,\,IU(x,\,y)\,\longrightarrow\,x\,\curlywedge\,y\,=\,U_{0}(x)\,=\,x\,\owedge\,y,\,x\,\curlyvee\,y\,=\,U_{0}(y)\,=\,x\,\ovee\,y)$, can be proved }
\item {$(IU(x,\,y,\,a,\,b),\,x\,\leq\,y,\,a,\,\leq\,b\,\longrightarrow\,x\,\curlywedge\,a\,\leq\,y\,\curlywedge\,b)$ can be verified in the same way as the above.}
\item {In the premise of  $(IU(x,\,y,\,a,\,b,\,x\,\ovee\,a,\,y\,\ovee\,b),\,x\,\leq\,y,\,a,\,\leq\,b\,\longrightarrow\,x\,\curlyvee\,a\,\leq\,y\,\curlyvee\,b)$, we require $IU(x\,\ovee\,a)$ and $IU(y\,\ovee\,b)$. So if $A\,\in\,x\,\curlyvee\,a$ and $B\,\in\,y\,vee\,b$, then if $C\,\in\, U_{0}([A])$ and $E\,\in\,U_{0}([B])$, then it can be checked that $C\,\subseteq\,E$. The proof consists in continuing the verification for the other operators. }
\item {For the $\Rightarrow$ part of  $t(x)\,=\,x$ if and only if $\neg(IU(x)),\,s(\breve{U}(x))\,=\,0,\,L_{0}(x)\,\lneq\,x,\,(\forall{y})(y\,\lneq\,x\,\longrightarrow\,ty\,=\,0)$, $(0\,\lneq\,a,\,b,\,c\,\lneq\,x\,\longrightarrow\,a\,=\,b\; \mathrm{or} \;b\,=\,c\; \mathrm{or} \;c\,=\,a)$ note that $x$ is the class of a two element set that is not contained in a single block (this is the last sentence by definition). Using those two elements as representatives of the class, it can be seen that its $\breve{U}$ approximation cannot be a disjoint union of blocks. So $s(\breve{U}(x))\,=\,0$ and $L_{0}(x)\,\lneq\,x$.

For the converse, a contradiction argument using the class of a one-element, and \emph{more than two element} sets yields the result.}
\item { If $A\,\in\,x$, then $\sim \minuso [A]\,=\,[S\,\setminus\,(S\,\setminus\,A^{u0})^{\breve{u}}]$ (assuming that $IU(x)$ holds), while $\minuso \sim [A]\,=\,[S\,\setminus\,(S\,\setminus\,A^{\breve{u}})^{u0}]$. Observe that $S\,\setminus\,(S\,\setminus\,A^{u0})^{\breve{u}}\,\subseteq\,S\,\setminus\,(S\,\setminus\,A^{\breve{u}})^{u0}$ holds in particular for $A$. As the operations are well-defined, we have
$(IU(x)\,\longrightarrow\,\sim \minuso x\,\leq\,\minuso \sim x)$.}
\item { If $A\,\in\,x$, then $\sim [A]\,=\,[S\,\setminus\,A^{\breve{u}}]$. Again $\sim L_{0}([A])\,=\,[S\,\setminus\,(A^{l0})^{\breve{u}}]$. As $S\,\setminus\,A^{\breve{u}}\,\subseteq\,S\,\setminus\,(A^{l0})^{\breve{u}}$, so  $\sim x \,\leq\,\sim L_{0}(x)$ ; 

$\sim 0\,=\,[S\,\setminus\,\emptyset^{\breve{u}}]\,=\,[S]\,=\,1$ ; Similarly $\sim 1 \,=\,0$}
\item {The proof of the two statements $x\,\leq\,\sim \sim x$ and $(IU(x)\,\longrightarrow\,\sim U_{0} (x)\,\leq\,\sim x)$ is similar to that of the above.}
\item {Let $A\,\in\,x$, then $\sim \breve{U}[A]\,=\,[S\,\setminus\,(A^{\breve{u}})^{\breve{u}}]$ and $S\,\setminus\,(A^{\breve{u}})^{\breve{u}}\,\subseteq\,S\,\setminus\,A^{\breve{u}}$, so $\sim \breve{U}(x)\,\leq\,\sim x $ ; The second part $(IU(x)\,\longrightarrow\,\sim \breve{U} (x)\,\leq\, \sim U_{0}(x)$ follows from the first part provided $IU(x)$ holds.}
\item {The proof of $\sim x \,\leq\,\sim \breve{L}(x)\,\leq\,\sim L_{0}(x)$ is similar to the proof of the above. 

For the second part, if $A\,\in\,x$, then $IN(x)$ will ensure that $\minuso [A]$ is defined and $\minuso [A]\,=\,[S\,\setminus\,A^{u0}]$ and $S\,\setminus\,A^{u0}\,\subseteq\,S\,\setminus\,(A^{l0})^{u0}$. So  $(IN(x)\,\longrightarrow\,\minuso x\,\leq\,\minuso L_{0}(x))$}
\item {$(IU(\sim x),\,IU(x)\,\longrightarrow\,\sim \minuso x \,\leq\, \minuso \sim x)$ is not very hard to prove. 

$\minuso 0\,=[S\,\setminus\,\emptyset]\,=\,\,1$. 

$\neg IN(1)$ in general holds because the 0-upper approximation of $S$ will not be defined unless the tolerance relation $T$ is an equivalence in the first place.}
\item {The proof of $(IU(x)\,\longrightarrow\,\minuso x\,\leq\,\minuso L_{0}(x))$ and\\ $(IN(x)\,\longrightarrow\,\minuso \breve{U}(x)\,\leq\,\minuso x,\,\minuso U_{0}(x)\,=\,\minuso x)$ is fairly direct}
\item {Note that $s(x)\,=\,x$ is the same thing as saying that $x$ is a union of disjoint blocks. So the statement, $\forall{y}(y\,\lneq\,x\,\longrightarrow\,L_{0}(y)\,\lneq\,L_{0}(x))$ and $L_{0}(x)\,=\,U_{0}(x)\,=\,x$ if and only if $s(x)\,=\,x$, else $s(x)\,=\,0$ holds.}
\end{enumerate}
The proof of the last two statements is easy. \qed
\end{proof}

\begin{proposition}[Implication-Like Operations]
In $\maltese (S)$, if we define $x\,\rightsquigarrow\,y\,=(\sim x)\,\sqcup\,(\breve{U} y)$ and $x\,\rightarrowtail\,y\,=\,(\minuso x)\sqcup\,(U_{0} y)$ if defined, then:
\begin{enumerate}[(a)]
\item {$(IU(x)\,\longrightarrow\,x\,\rightarrowtail\,x\,=\,1)$;  $(IU(x)\,\longrightarrow\,U_{0}(x)\,\leq\,1\,\rightarrowtail\,x)$} 
\item {$(IU(x,\,y)\,\longrightarrow\,U_{0}(x)\,\leq\,x\,\rightarrowtail\,(y\,\rightarrowtail\,x) )$;  $(IU(x)\,\longrightarrow\,x\,\rightarrowtail\,0\,=\,\minuso(x)) $}
\item {$x\,\rightsquigarrow\,x\,=\,1 $;  $x\,\,\leq\,(1,\rightsquigarrow\,x) $; $x\,\rightsquigarrow\,0\,=\,\sim\,x $}
\item {  $(x\,\rightsquigarrow\,y)\,\sqcup\,(x\,\rightsquigarrow\,z)\,\leq\,(x\,\rightsquigarrow\,(y\,\sqcup\,z)) $} 
\item {$(x\,\rightsquigarrow\,z)\,\sqcap\,(y\,\rightsquigarrow\,z)\,\leq\,(x\,\sqcup\,y)\,\rightsquigarrow\,z $ \qed} 
\end{enumerate}
\end{proposition}
\section{Representation Theorems}
\begin{definition}
By an \emph{abstract essential $\lambda$-rough partial algebraic system} (AER) we will mean a partial algebraic system of the form \[S\,=\,\left\langle \underline{S},\,\leq ,\,\sqcup,\,\owedge,\,L_{0},\,U_{0},\,\breve{L},\,\breve{U},\,\sim,\,\minuso,\,t,\,0,\,1,\,(2,\,2,\,2,\,1,\,1,\,1,\,1,\,1,\,1,\,1,\,0,\,0) \right\rangle \] that satisfies all of the following (we assume that the operation $\sqcup$ is complete and the derived operations $\ovee,\,\curlyvee,\,\curlywedge,\,s $ and derived predicates $\lneq,\,IU,\,IN $ are defined via:
\begin{itemize}
\item {If the RHS is defined then and only then $x\,\ovee\,y\,{=}\,U_{0}(x)\,\sqcup\,U_{0}(y)$}
\item {$x\,\leq\,y,\,\neg(x\,=\,y)$ if and only if $x\,\lneq\,y$}
\item {$x\,\curlyvee\,y\,\stackrel{w}{=}\,U_{0}(x\,\ovee\,y)$ ; $x\,\curlywedge\,y\,\stackrel{w}{=}\,L_{0}(x\,\owedge\,y)$ ; $IN(x)$ if and only if $\minuso x\,=\,\minuso x$.}
\item {$\forall{y}(y\,\lneq\,x\,\longrightarrow\,L_{0}(y)\,\lneq\,L_{0}(x))$ and $L_{0}(x)\,=\,U_{0}(x)\,=\,x$ if and only if $s(x)\,=\,x$, else $s(x)\,=\,0$}
\item {$IU(x)$ if and only if $U_{0}(x)\,=\,U_{0}(x)$. Further we will write $IU(a, b,..)$ for $IU(a),\, IU(b),\ldots$}
\end{itemize}
\begin{enumerate}
\item {$x\,\curlyvee\,y\,\stackrel{w}{=}\,y\,\curlyvee\,x$ ; $x\,\curlywedge\,y\,\stackrel{w}{=}\,y\,\curlywedge\,x$ ; $x\,\curlyvee\,x\,\stackrel{w}{=}\,U_{0} U_{0}(x)$ ; $x\,\curlywedge\,x\,\stackrel{w}{=}\,U_{0}(x)$}
\item {$L_{0}(x)\,\leq\,\breve{L}(x)\,\leq\,x\,\breve{U}(x)$ ; $(IU(x)\,\longrightarrow\,x\,\leq\,U_{0}(x)\,\leq\,U_{0} U_{0}(x))$}
\item {$L_{0}L_{0}(x)\,=\,L_{0}(x)$ ; $\breve{L}L_{0}(x)\,=\,L_{0}(x)$ ; $L_{0}\breve{L}(x)\,\leq\,\breve{L}(x)$}
\item {$\breve{L}\breve{L}(x)\,=\,\breve{L}(x)$ ; $(IU(x)\,\longrightarrow\, L_{0}U_{0}(x)\,=\,U_{0}(x))$ ; $L_{0}(x)\,\leq\,U_{0}L_{0}(x)$}
\item {$\breve{U}(x)\,\leq\,\breve{U}\breve{U}(x)$ ; $(IU(x)\,\longrightarrow\,x\,\leq\,U_{0}(x)\,\leq\,\breve{U}(x)\,\leq\,\breve{U}U_{0}(x)\,\leq\,\breve{U}\breve{U}(x))$}
\item {$(x\,\leq\,y,\,IU(x,\,y)\,\longrightarrow\,x\,\curlywedge\,y\,=\,U_{0}(x)\,=\,x\,\owedge\,y,\,x\,\curlyvee\,y\,=\,U_{0}(y)\,=\,x\,\ovee\,y)$}
\item {$(IU(x,\,y,\,a,\,b),\,x\,\leq\,y,\,a,\,\leq\,b\,\longrightarrow\,x\,\curlywedge\,a\,\leq\,y\,\curlywedge\,b)$}
\item {$(IU(x,\,y,\,a,\,b,\,x\,\ovee\,a,\,y\,\ovee\,b),\,x\,\leq\,y,\,a,\,\leq\,b\,\longrightarrow\,x\,\curlyvee\,a\,\leq\,y\,\curlyvee\,b)$}
\item {$(x\,\leq\,y\,\longrightarrow\,L_{0}(x)\,\leq\,L_{0}(y),\,\breve{U}(x)\,\leq\,\breve{U}(y),\,\breve{L}(x)\,\leq\,\breve{L}(y))$}
\item {$(x\,\leq\,y,\,IU(x)\,\longrightarrow\,U_{0}(x)\,\leq\,U_{0}(y))$}
\item {$t(x)\,=\,x$ if and only if $\neg(IU(x)),\,s(\breve{U}(x))\,=\,0,\,L_{0}(x)\,\lneq\,x,\,(\forall{y})(y\,\lneq\,x\,\longrightarrow\,ty\,=\,0)$ and $(0\,\lneq\,a,\,b,\,c\,\lneq\,x\,\longrightarrow\,a\,=\,b\; \mathrm{or} \;b\,=\,c\; \mathrm{or} \;c\,=\,a)$}
\item {$(IU(x)\,\longrightarrow\,\sim \minuso x\,\leq\,\minuso \sim x)$ ; $\sim x \,\leq\,\sim L_{0}(x)$ ; $\sim 0\,=\,1$ ; $\sim 1 \,=\,0$}
\item {$(IU(\sim \sim x)\,\longrightarrow\,x\,\leq\,\sim \sim x)$ ; $(IU(x)\,\longrightarrow\,\sim U_{0} (x)\,\leq\,\sim x)$}
\item {$\sim \breve{U}(x)\,\leq\,\sim x $ ; $(IU(x)\,\longrightarrow\,\sim \breve{U} (x)\,\leq\, \sim U_{0}(x)$}
\item {$\sim x \,\leq\,\sim \breve{L}(x)\,\leq\,\sim L_{0}(x)$ ; $(IN(x)\,\longrightarrow\,\minuso x\,\leq\,\minuso L_{0}(x))$}
\item {$(IU(\sim x),\,IU(x)\,\longrightarrow\,\sim \minuso x \,\leq\, \minuso \sim x)$ ; $\minuso 0\,=\,1$ ; $\neg IN(1)$}
\item {$(IU(x)\,\longrightarrow\,\minuso x\,\leq\,\minuso L_{0}(x))$ ; $(IN(x)\,\longrightarrow\,\minuso \breve{U}(x)\,\leq\,\minuso x,\,\minuso U_{0}(x)\,=\,\minuso x)$}
\item {$\forall{y}((0\,\leq\,y\,\leq\,x\,\longrightarrow\,y\,=\,0\; \mathrm{or} \;y\,=\,x),\,\longrightarrow\,\bigvee (s(z)\,=\,z,\,x\,\leq\,z)) $, where $\bigvee$ indicates disjunction over the entire set $S$}
\item {$(s(x)\,=\,x,\,x\,\lneq\,y\,\longrightarrow\,s(y)\,=\,0)$ ; $(s(x)\,=\,x,\,y\,\lneq\,x\,\longrightarrow\,s(y)\,=\,0)$}
\item {$(\forall{x})\,0\,\leq\,x\,\leq\,1$ ; $x\,\sqcup\,y\,=\,y\,\sqcup\,x$ ; $x\,\sqcup\,x\,=\,x$ ; }
\item {$(x\,\leq\,y\,\longrightarrow\,x,\,y\,\leq\,x\,\sqcup\,y)$ ; $\breve(U)(x)\,\sqcup\,\breve{U}(y)\,\leq\,\breve{U}(x\,\sqcup\,y)$}
\end{enumerate}
\end{definition}

\begin{theorem}
Given an abstract essential $\lambda$-rough partial algebra $S$ there exists a tolerance approximation space and a choice perspective that ensures that its algebraic semantics is isomorphic to $S$. 
\end{theorem}
\begin{proof}
Our abridged proof has three components (roughly). The first concerns the reconstruction of the tolerance approximation space, the second part of the choice perspective and the third part of compatibility builds into the first two.

Let $S$ be an AER as in the above definition. Then the statements $11,\,16-21$ and the definitions of $\lneq,\,s,\,t$, ensure that we can reconstruct a tolerance $T$ on a set $K$ (corresponds to $1$) by the representation theorem of tolerance relations (see \cite{CZ}). $\sqcup$ is needed for getting the set $K$ in a easier way. We do not have a full representation here.

Both the operations $L_{0}$ and $U_{0}$ permit the isolation of the choice function used as blocks can be identified through the function $s$ and combined via $\sqcup$, while maximal unions of mutually disjoint blocks can be identified and any union of blocks is constructible. \qed
\end{proof}

\section{Further Directions: Modal Connections}

A relational structure, in particular a TAS, can be associated with many modal logics by way of induced global operations through suitable modification of the theory of Tarski algebras (see \cite{CSA}). Importantly, the semantics developed above can be combined with such modal semantics (\cite{AM99}). Choice functions \emph{induced by modalities} are also considered in the same paper by the present author. A stronger abstract representation theorem is proved for many of the contexts; This is partly due to the partial representation results in \cite{CSA}. See the same for undefined concepts in this section. 

In particular, if $\mathbb{S}$ is the set of blocks of a TAS $S$, then $K\,=\,\left\langle \underline{S},\,\mathbb{S}\right\rangle $ is a \emph{dense Tarski set}.  Let $\Delta(S)$ be the set of subsets that are unions of 'a complement of a block and a subset of it'. Then $X\,=\,\left\langle \Delta(S),\,\Rightarrow,\,\underline{S} \right\rangle$ is a \emph{Tarski Subalgebra} of $\wp (S)$, where $\Rightarrow$ is defined via $U\,\Rightarrow\,V\,=\,(\underline{S}\,\setminus\,U)\,\cup\,V$. The set of maximal filters of $X$ allow a representation theorem in the finite case. The form of $\Delta(S)$ in the above is actually dictated by the choice of $\mathbb{S}$, which may be taken to be the collection of unions of disjoint blocks among other things.  

\begin{proposition}
If the rough equivalence $\sigma$ of the previous section is applied on $\Delta(S)$, then the resulting classes are of two types: those that correspond to complements of a block and those that correspond to unions of complement of a block and a little more. \qed
\end{proposition}

The implication-like operations on an AER are clearly different from that of modal Tarski algebras and the best way of getting to a unified semantics is still not fully solved (even in \cite{AM99}). We can also adjoin arbitrary modal operations to an AER. The connections with other similarity semantics and logics can be found in more detail in the same paper.   

\bibliographystyle{splncs}
\bibliography{newsem33.bib}
\end{document}